\documentclass[a4paper,11pt]{amsart}
\usepackage{amssymb}
\usepackage{mathrsfs}
\usepackage{amsmath,amssymb,amsthm,latexsym,amscd,mathrsfs}
\usepackage{indentfirst}
\numberwithin{equation}{section} \theoremstyle{plain}
\newtheorem{thm}{Theorem}[section]

\newtheorem{prop}{Proposition}[section]

\newtheorem{cor}{Corollary}[section]
\theoremstyle{definition}

\theoremstyle{remark}
\newtheorem{rem}{Remark}[section]

\newtheorem{ack}{Acknowledgements}   

\title[A note on Ribaucour transformations]
{A note on Ribaucour transformations in Lie sphere geometry}
\author[J.Q. Ge]{Jianquan Ge}
\address{School of Mathematical Sciences, Laboratory of Mathematics and Complex Systems, Beijing Normal
University, Beijing 100875, P.R. CHINA}
\address{Current address: Mathematisches Institut der Universit\"{a}t zu K\"{o}ln - Weyertal 86-90
- 50931 K\"{o}ln, Germany}
\email{jqge@bnu.edu.cn}
\thanks {The author is partially supported by the NSFC (No.11001016), the SRFDP (No. 20100003120003), and a research fellowship from the Alexander von Humboldt Foundation.}

\date{}
\begin{document}
\maketitle

\begin{abstract}
Following Burstall and Hertrich-Jeromin we study the Ribaucour transformation of Legendre submanifolds in Lie sphere geometry. We give an explicit parametrization of the resulted Legendre submanifold $\hat{F}$ of a Ribaucour transformation, via a single real function $\tau$ which represents the regular Ribaucour sphere congruence $s$ enveloped by the original Legendre submanifold $F$.
\end{abstract}

\section{Introduction}
Classically the Ribaucour transformation and its generalizations in submanifold geometry have been intensively studied (cf. \cite{BH, DT1, DT2} and references therein).

Geometrically a \emph{Ribaucour transform} of an immersion $f:\Sigma^n\rightarrow N^{m+1}$ in a space-form gives another immersion $\hat{f}:\Sigma^n\rightarrow N^{m+1}$ with the properties: (1) for each $p\in\Sigma^n$ there is an $n$-sphere in $N$ tangent to both $df(T_p\Sigma)$ and $d\hat{f}(T_p\Sigma)$; (2) the shape operators of $f$ and $\hat{f}$ commute. Although it has a clear geometric background, an intuitive picture or parametrization for a ``non-trivial" (for example: not a parallel transformation, $f$ is not umbilical, \emph{etc.}) Ribaucour transform is rather subtle and involved, even in the classical surface case. Via a beautiful correspondence to commuting Codazzi tensors and Combescure transforms, Dajczer and Tojeiro \cite{DT2} were able to give a parametrization for general Ribaucour transforms of a simply connected submanifold in $\mathbb{R}^{n+p}_s$ (indefinite space form of index $s$). They also successfully generalized the Bianchi Permutability Theorem to higher dimensional and co-dimensional case. Due to the Lie invariant nature and complex PDEs, it seems still cumbersome to imagine how a Ribaucour transform act in submanifold geometry.

An elegant treatment of Ribaucour transformations in its natural context: Lie sphere geometry has been initiated by Burstall and Hertrich-Jeromin \cite{BH}, which contains also a simple conceptual proof of the generalized Bianchi Permutability Theorem, which states that, given two Ribaucour transforms with certain condition satisfied, there exist two $\mathbb{R}P^1$-families, the \emph{Demoulin families} of Legendre maps such that each map in one family is a Ribaucour transform of each map in the other. They called a sphere congruence $s$ enveloped by two legendre maps $F,\hat{F}$ a \emph{Ribaucour sphere congruence} if the quotient bundle $(F+\hat{F})/s$ is flat with respect to the natural metric connection induced from the exterior differential $d$ (see for details later). Then in this case, they called $F,\hat{F}$ \emph{Ribaucour transforms} of each other and a \emph{Ribaucour pair}. A direct parametrization in Lie sphere geometry as that \cite{DT2} did in submanifold geometry mentioned above was left unknown, though the definitions and results were proved to descend to those in submanifold geometry.

It is our aim in this paper to present explicitly the parametrization of the Ribaucour sphere congruence and the Ribaucour transformation from a given Legendre submanifold, under a relatively mild and natural condition of regularity of the sphere congruence. Besides the given Legendre map $F$, we use only a real function $\tau$ on the manifold to determine the parametrization both for the Ribaucour sphere congruence $s$ and the new Legendre submanifold $\hat{F}$ obtained by the Ribaucour transformation. In particular, the Ribaucour condition can be identified with the closedness of a $1$-form determined by $\tau$. Moreover, we can recover $F$ from $\hat{F}$ and $s$ (or the same representative function $\tau$) by the same process. As an application, we can give explicit parametrization for the maps of the Demoulin families in the Bianchi Permutability Theorem in terms of the representative functions $\tau_0,\tau_1$ of the given two Ribaucour transformations.
\section{Regular sphere congruences of Legendre submanifolds}
Denote by $\mathbb{R}^{m+2,2}$ the $(m+4)$-dimensional vector space with metric $(,)$ of signature $(m+2,2)$. Let $\mathcal{Q}$ be the Lie quadric, \emph{i.e.}, the projective light-cone of $\mathbb{R}^{m+2,2}$, and $\mathcal{Z}\cong T_1S^{m+1}$ the contact manifold of projective lines on $\mathcal{Q}$. A map from a manifold $M$ to $\mathcal{Q}$ (resp. $\mathcal{Z}$) can be also recognized as a line (resp. $2$-plane) subbundle of the product $M\times\mathbb{R}^{m+2,2}$ (cf. \cite{BH,Ce}).

Recall \cite{BH} that a \emph{sphere congruence} is a map $s: M^m\rightarrow \mathcal{Q}$ of an $m-$manifold $M^m$, a Legendre map is a map $F: M^m\rightarrow\mathcal{Z}$ satisfying the contact condition
\begin{equation}\label{contact}
(d\sigma_0,\sigma_1)\equiv0 \quad \emph{for all sections $\sigma_0,\sigma_1$ of $F$,}
\end{equation}
and a Legendre map $F$ \emph{envelops} a sphere congruence $s$ if $s(p)\subset F(p)$ for all $p\in M^m$. It is well known \cite{Ce} that a Legendre map $F$ corresponds to its \emph{point sphere map} $F^0:M^m\rightarrow\mathcal{Q}$ and \emph{great sphere map} $F^1:M^m\rightarrow\mathcal{Q}$ by $F=[F^0,F^1]$, the projective line through $F^0:=\langle f+t_0\rangle$, $F^1:=\langle\xi+t_1\rangle$ on $\mathcal{Q}$, where $t_0,t_1$ are fixed orthogonal unit time-like vectors, $f: M^m\rightarrow S^{m+1}\subset \mathbb{R}^{m+2}=\langle t_0, t_1\rangle^{\bot}$ and $\xi:M^m\rightarrow S^{m+1}$ are the \emph{spherical projection} and the \emph{spherical field of unit normals} of $F$ respectively. Moreover, they satisfy
\begin{equation}\label{contactpt}
(f,\xi)=0,\quad (df,\xi)=(f,d\xi)=0
\end{equation}
since for $F^0$ and $F^1$ are in oriented contact and the contact condition (\ref{contact}) of $F$. The Legendre map $F$ is called a \emph{Legendre submanifold} if it is further an immersion, or equivalently, $df(X),d\xi(X)$ do not simultaneously vanish for all nonzero vector $X\in TM$.

Given a Legendre submanifold $F:M^m\rightarrow \mathcal{Z}$, we can express each its enveloping sphere congruences $s$ which is pointwise distinct with the point sphere map $F^0$ by
\begin{equation}\label{spherecong}
s=\langle\xi-\tau f-\tau t_0+t_1\rangle=:\langle\sigma\rangle\subset F^1-\tau F^0\subset F \quad \emph{for $\tau\in C^{\infty}(M)$.}
\end{equation}
Then $s$ is called a \emph{regular} sphere congruence enveloped by $F$ if $-d\xi+\tau df$ is non-degenerate everywhere. When $F$ is induced from a hypersurface (\emph{i.e.}, $f$ is an immersion), this is precisely the condition that the sphere congruence contains no curvature spheres, or equivalently, $-\tau$ differs with principal curvatures of the hypersurface $f$ in direction $\xi$. Now since $df, d\xi$ do not vanish simultaneously at any point, regular sphere congruences enveloped by $F$ indeed exist globally. Next we will show that a regular sphere congruence $s$ of a Legendre submanifold $F$ is naturally enveloped by another Legendre submanifold $\hat{F}=[\hat{F}^0,\hat{F}^1]$ by giving its spherical projection $\hat{f}$ and spherical field of unit normals $\hat{\xi}$.

First we split $\hat{f}:M^m\rightarrow S^{m+1}$ as
\begin{equation}\label{fhat}
\hat{f}=af+b\xi+(1-a)\check{f},
\end{equation}
where $a,b\in C^{\infty}(M)$ and $\check{f}\in \langle f,\xi\rangle^{\bot}\subset\mathbb{R}^{m+2}$ are to be determined. Since $s$ is required to be enveloped by $\hat{F}$, \emph{i.e.}, $s\subset \hat{F}=[\hat{f}+t_0,\hat{\xi}+t_1]$, $\hat{\xi}:M^m\rightarrow S^{m+1}$ has already been determined as
\begin{equation}\label{xihat}
\hat{\xi}=\xi-\tau f+\tau \hat{f}.
\end{equation}
By the definitions (\ref{fhat}, \ref{xihat}) and the contact condition (\ref{contactpt}), we have the following required equations:
\begin{equation}\label{contactpt-fhat}
\begin{array}{llll}
1=(\hat{f},\hat{f})=a^2+b^2+(1-a)^2(\check{f},\check{f}),&&& \\
0=(\hat{f},\hat{\xi})=b-a\tau+\tau,&&&\\
1=(\hat{\xi},\hat{\xi})=1+2\tau^2-2a\tau^2+2b\tau,&&&\\
0=(d\hat{f},\hat{\xi})=db-\tau da+(1-a)(-d\xi+\tau df, \check{f}).&&&
\end{array}
\end{equation}
Now we introduce a Riemannian metric $\langle,\rangle_{-}$ on $M^m$ by
\begin{equation}\label{metricbar}
\langle X,Y\rangle_{-}:=((-d\xi+\tau df)(X), (-d\xi+\tau df)(Y)) \quad \emph{for $X,Y\in TM$,}
\end{equation}
and denote the corresponding Levi-Civita connection (covariant derivative, gradient) by $\bar{\nabla}$. This is well-defined by regularity assumption. Then easy computation of (\ref{contactpt-fhat}) shows
\begin{equation}\label{checkf}
\begin{array}{ll}
a=1-\frac{2}{\tau^2+\mu^2+1}, & b=\tau(a-1)=\frac{-2\tau}{\tau^2+\mu^2+1},\\
\check{f}=(-d\xi+\tau df)(\bar{\nabla}\tau), & \mu^2:=(\check{f},\check{f})=\langle\bar{\nabla}\tau,\bar{\nabla}\tau\rangle_{-}.
\end{array}
\end{equation}
Finally, taking differential of (\ref{xihat}) shows
\begin{equation}\label{pair}
-d\hat{\xi}+\tau d\hat{f}=-d\xi+\tau df +(f-\hat{f})d\tau
\end{equation}
is non-degenerate everywhere whence $\hat{F}$ is a Legendre submanifold and s is a regular sphere congruence enveloped both by $F$ and $\hat{F}$.

Conversely, we can reconstruct $F$ from $\hat{F}$ in the same way with the same functions $\tau$ and $\mu$ (and hence also $a,b$) due to (\ref{xihat}, \ref{pair}). In fact, direct calculations show that the metric induced by $-d\hat{\xi}+\tau d\hat{f}$ is the same as that by $-d\xi+\tau df$ in (\ref{metricbar}) and hence the map $\check{\hat{f}}$ in the decomposition of $f$ with respect to $\hat{f}, \hat{\xi}$ as in (\ref{fhat}) can be written as
\begin{equation}\label{checkhatf}
\check{\hat{f}}=\check{f}+\mu^2(f-\hat{f})=\check{f}+\mu^2(1-a)(f+\tau\xi-\check{f})
\end{equation}
from which it follows that $\check{\hat{f}}$ has the same length $\mu$ as $\check{f}$.

 We conclude these as
\begin{prop}\label{reg-cong-pair}
Given a Legendre submanifold $F=[f+t_0,\xi+t_1]:M^m\rightarrow\mathcal{Z}$ and a regular sphere congruence $s$ represented by a function $\tau\in C^{\infty}(M)$ as in $(\ref{spherecong})$, there exists another unique Legendre submanifold $\hat{F}=[\hat{f}+t_0,\hat{\xi}+t_1]:M^m\rightarrow\mathcal{Z}$ as parametrized in $(\ref{fhat},\ref{xihat},\ref{checkf})$ enveloping $s$, \emph{i.e.}, $s=F\cap\hat{F}$. Conversely we can reconstruct $F$ from $(\hat{F}, s)$ in the same way with the same functions $\tau,\mu,a,b$. In conclusion, any two of $\{F,s,\hat{F}\}$ determine the third.
\end{prop}
\begin{rem}
When $f:M^m\rightarrow S^{m+1}$ is a hypersurface with $\xi$ its oriented unit normal vector field, we can use the induced metric $\langle,\rangle_f$ instead of the metric defined in (\ref{metricbar}) so that $\check{f}$ in (\ref{checkf}) can be rewritten more familiarly as
\begin{equation}\label{checkf-hyp}
\check{f}=df\circ (A^{\xi}+\tau Id)^{-1}(\nabla^f \tau),
\end{equation}
where $A^{\xi}:TM\rightarrow TM$ is the shape operator of $f$ and $\nabla^f \tau$ is the gradient with respect to the induced metric of $f$.
\end{rem}
\begin{rem}\label{nonregular}
When $s$ is not regular, \emph{i.e.}, $-d\xi+\tau df$ is degenerate somewhere, the last equation of (\ref{contactpt-fhat}) for $\check{f}$ may have no solution or have more than one solution. In the latter case, the images of these new Legendre maps enveloping $s$ may still coincide with each other, since, roughly speaking, these maps would only differ with each other as flows in principal directions.
\end{rem}

\section{Ribaucour condition}
Let $\{F,s,\hat{F}\}$ be the triple of two Legendre submanifolds enveloping a regular sphere congruence with representative function $\tau\in C^{\infty}(M)$ as in Proposition \ref{reg-cong-pair}. In this section we will give a characterization of the Ribaucour condition by the closedness of a $1$-form on $M$ determined by $\tau$ starting from either of $\{F, \hat{F}\}$.

First we recall the original definition of the Ribaucour transformation in Lie sphere geometry given in \cite{BH}. Let $\mathcal{N}_{F,\hat{F}}:=(F+\hat{F})/s$ be the quotient bundle of rank $2$ over $M$. By the contact condition (\ref{contact}), $\mathcal{N}_{F,\hat{F}}$ is a subbundle of $s^{\bot}/s$ whence it inherits a metric of signature $(1,1)$ and a metric connection $\nabla^{F,\hat{F}}(\nu+s):=\pi(d\nu+s)$ where $\pi:s^{\bot}/s\rightarrow\mathcal{N}_{F,\hat{F}}$ is the orthogonal projection. As mentioned in the Introduction, $s$ (resp. $F,\hat{F}$) is called a \emph{Ribaucour sphere congruence} (resp. \emph{Ribaucour transforms} of each other or \emph{Ribaucour pair}) if  $\nabla^{F,\hat{F}}$ is flat. Meanwhile, there is an alternative definition presented as follows (see also \cite{BC}).
Consider the point sphere maps $F^0=\langle f+t_0\rangle, \hat{F}^0=\langle\hat{f}+t_0\rangle$ of $F,\hat{F}$. Then $\langle F^0,\hat{F}^0\rangle$, the space spanned by $F^0,\hat{F}^0$ at each point, is a rank $2$ subbundle with metric and connection $\nabla^{F^0,\hat{F}^0}$ induced from $M^m\times \mathbb{R}^{m+2,2}$. It is easily seen that $\mathcal{N}_{F,\hat{F}}\cong\langle F^0,\hat{F}^0\rangle$ by a metric, connection preserving isomorphism given by $$\nu+s\rightarrow \nu+(\nu,t_1)\sigma,\quad (\emph{here $\langle\sigma\rangle=s$ as in $(\ref{spherecong})$})$$
which indicates $\mathcal{N}_{F,\hat{F}}$ is flat if and only if $\langle F^0,\hat{F}^0\rangle$ is flat. As shown in \cite{BH}, the latter condition is easy to characterize by the vanishing of the curvature form, which, by the Gauss equation and the flatness of $d$, is equal to
\begin{equation}\label{curvature}
(\beta(f+t_0)\wedge\beta(\hat{f}+t_0))=0,
\end{equation}
where $\beta\in\Omega^1_{M^m}\otimes Hom(\langle F^0,\hat{F}^0\rangle, \langle F^0,\hat{F}^0\rangle^{\bot})$ is the second fundamental form of $\langle F^0,\hat{F}^0\rangle$ given by
\begin{equation}\label{beta}
\beta \psi=d\psi-\nabla^{F^0,\hat{F}^0}\psi \quad \emph{for $\psi\in\Gamma\langle F^0,\hat{F}^0\rangle$,}
\end{equation}
$``\wedge"$ acts on the $1$-form-factors and $(\cdot\wedge\cdot)$ takes inner product on the vector-factors.

Now we define $1$-forms $\alpha_{\tau},\hat{\alpha}_{\tau}\in\Omega^1_{M^m}$ for $\tau\in C^{\infty}(M)$ with respect to $F,\hat{F}$, respectively, by (\ref{checkf}, \ref{checkhatf}) as the following:
\begin{equation}\label{alpha}
\begin{array}{ll}
\alpha_{\tau}:=(df, (d\xi-\tau df)(\bar{\nabla}\tau))=(df,-\check{f})=(a-1)^{-1}(df,\hat{f}),&\\
\hat{\alpha}_{\tau}:=(d\hat{f}, (d\hat{\xi}-\tau d\hat{f})(\bar{\nabla}\tau))=(d\hat{f},-\check{\hat{f}})=(a-1)^{-1}(d\hat{f},f).
\end{array}
\end{equation}
Note that since $(f,\hat{f})=a$, we have
\begin{equation}\label{alphasum}
\alpha_{\tau}+\hat{\alpha}_{\tau}=d\ln(1-a).
\end{equation}
Direct calculations show that
\begin{equation}\label{connectionN}
\nabla^{F^0,\hat{F}^0}(f+t_0)=(f+t_0)\alpha_{\tau}, \quad \nabla^{F^0,\hat{F}^0}(\hat{f}+t_0)=(\hat{f}+t_0)\hat{\alpha}_{\tau},
\end{equation}
from which, applying formulae (\ref{beta}, \ref{alpha}, \ref{alphasum}), we get
\begin{equation*}
(\beta(f+t_0)\wedge\beta(\hat{f}+t_0))=(1-a)d\alpha_{\tau}.
\end{equation*}
Hence, from (\ref{curvature}) we arrive at the equivalent Ribaucour condition: $d\alpha_{\tau}=0$, or equivalently, $d\hat{\alpha}_{\tau}=0$ by (\ref{alphasum}).

Combining the construction in last section, we obtain
\begin{thm}\label{Ricond}
Given a Legendre submanifold $F=[f+t_0,\xi+t_1]:M^m\rightarrow\mathcal{Z}$ and a function $\tau\in C^{\infty}(M)$ satisfying $-d\xi+\tau df$ is non-degenerate everywhere,
let $s:M^m\rightarrow \mathcal{Q}$ be the map defined in $(\ref{spherecong})$
and $\hat{F}=[\hat{f}+t_0,\hat{\xi}+t_1]:M^m\rightarrow \mathcal{Z}$ be the map defined by $(\ref{fhat},\ref{xihat},\ref{checkf})$.
 Then $s$ (resp. $\hat{F}$) is a regular Ribaucour sphere congruence (resp. Ribaucour transform) of $F$ if and only if $\alpha_{\tau}$ defined in $(\ref{alpha})$ is a closed $1$-form on $M$.
 Moreover, the same conclusions hold if we interchange $\hat{F}$ and $F$ with the same $\tau, s$.
\end{thm}
In fact, we have shown a correspondence between the subset of functions $\tau\in C^{\infty}(M)$ satisfying (1) $-d\xi+\tau df$ is non-degenerate; (2) $d\alpha_{\tau}=0$,
and the set of regular Ribaucour sphere congruences (``regular" Ribaucour transforms) of a given Legendre submanifold  $F=[f+t_0,\xi+t_1]:M^m\rightarrow\mathcal{Z}$.
Henceforth, we call such function $\tau$ a (regular) \emph{Ribaucour function} of $F$ and denote the set by $\mathcal{R}_F(M)$.
We also denote the regular Ribaucour sphere congruence $s$ (resp. Ribaucour transform $\hat{F}$) determined by $\tau$ as
$\mathcal{R}_{s_\tau}(F)$ (resp. $\mathcal{R}_{\tau}(F)$).
Then we have also shown
\begin{equation*}
\tau\in\mathcal{R}_F(M)\cap\mathcal{R}_{\hat{F}}(M), \quad \mathcal{R}_{s_\tau}(F)=\mathcal{R}_{s_\tau}(\hat{F}), \quad F=\mathcal{R}_{\tau}(\hat{F}).
\end{equation*}

As an application, we can present the Ribaucour functions $\{\tau_{\theta}\}\subset\mathcal{R}_F(M)$ for the Demoulin family $\{\hat{F}_{\theta}|e^{i\theta}\in\mathbb{R}P^1\}$ through the two given Ribaucour transforms $\hat{F}_0, \hat{F}_1$ of a Legendre submanifold $F$ in terms of the two pointwise distinct Ribaucour functions $\tau_0,\tau_1\in\mathcal{R}_F(M)$ corresponding to $\hat{F}_0, \hat{F}_1$. Before that, it is worthy to mention that by (\ref{curvature}, \ref{beta}, \ref{connectionN}), when $f$ is an immersion, there exist symmetric operators $r_{\tau_i}\in End(TM)$ (same as in \cite{BH}, symmetric with respect to $(df,df)$) such that $$d\hat{f_i}-(\hat{f_i}+t_0)\hat{\alpha}_{\tau_i}=(df-(f+t_0)\alpha_{\tau_i})\circ r_{\tau_i},$$
which derives a sufficient condition for the Bianchi Permutability Theorem as follows:
\begin{eqnarray}\label{Bianchicond}
&&[r_{\tau_0},r_{\tau_1}]=0,\quad \emph{or equivalently}\nonumber\\
&&((d\hat{f_0}-(\hat{f_0}+t_0)\hat{\alpha}_{\tau_0})\wedge(d\hat{f_1}-(\hat{f_1}+t_0)\hat{\alpha}_{\tau_1}))=0.
\end{eqnarray}
The latter also holds for general spherical projection $f$ by a similar argument as in \cite{BH}. This condition characterizes the flatness of the rank $4$ subbundle $V:=\hat{F}_0\oplus\hat{F}_1$ under the induced metric connection $\nabla$. Then the vector space of (locally) parallel sections of $V$ is isomorphic to $\mathbb{R}^{2,2}$ and its projective cone is a $(1,1)$-quadric $\mathcal{Q}^2$ isomorphic to $\mathbb{R}P^1\times\mathbb{R}P^1$. Note that by hypothesis, $\mathcal{N}_{F,\hat{F}_0}, \mathcal{N}_{F,\hat{F}_1}$ are flat so that we already have (locally) $\nabla$-parallel sections $\sigma_i\in\Gamma \mathcal{R}_{s_{\tau_i}}(F)\subset\Gamma F\subset\Gamma V$ looked as a $\nabla^{F,\hat{F}_j}$-parallel section in $\mathcal{N}_{F,\hat{F}_j}$, $\{i,j\}=\{0,1\}$ (cf.\cite{BH}). Then $\sigma_0,\sigma_1$ can be looked as two points in oriented contact on the quadric $\mathcal{Q}^2$ whence the projective line $[\sigma_0,\sigma_1]\cong\mathbb{R}P^1$ lies on $\mathcal{Q}^2$. In fact, the points of this line give exactly the Ribaucour sphere congruences $\mathcal{R}_{s_{\tau_{\theta}}}(F)$ of the Demoulin family $\{\hat{F}_{\theta}\}$ through $\hat{F}_0,\hat{F}_1$, since these Ribaucour sphere congruences, by the same argument as above, will have $\nabla$-parallel sections $\sigma_{\theta}\in\Gamma F$ in oriented contact with $\sigma_0,\sigma_1$ whence $\sigma_{\theta}\in[\sigma_0,\sigma_1]$. Therefore, we need only to figure out the $\nabla$-parallel sections $\sigma_0,\sigma_1$.

To give the Ribaucour functions globally, we also need to assume $M$ is simply connected as in \cite{BH}.
Now since $d\alpha_{\tau_i}=0$, there exist functions $\tilde{\tau}_i\in C^{\infty}(M)$ such that $\alpha_{\tau_i}=-d\tilde{\tau}_i$, $i=0,1.$
Write the sections $\sigma_0,\sigma_1$ as
$$\sigma_i:=u_i(\xi-\tau_i f-\tau_i t_0+t_1), \quad i=0,1,$$
where the functions $u_i$ are to be determined. It is not difficult to see that $\sigma_i$ is $\nabla$-parallel if and only if
$$(d\sigma_i, \hat{f}_j+t_0)=0, \quad \{i,j\}=\{0,1\}.$$
Then straightforward computation shows
$$u_i=\frac{C_ie^{\tilde{\tau}_j}}{\tau_i-\tau_j}, \quad \{i,j\}=\{0,1\}$$
for some nonzero constants $C_0,C_1$. Therefore, $\sigma_{\theta}\in\Gamma\mathcal{R}_{s_{\tau_{\theta}}}(F)$ corresponding to a linear combination of $\sigma_0,\sigma_1$, by (\ref{spherecong}), derives the family of Ribaucour functions as
\begin{equation}\label{Ribaufunctionfamily}
\tau_{\theta}=\frac{\cos\theta e^{\tilde{\tau}_1}\tau_0+\sin\theta e^{\tilde{\tau}_0}\tau_1}{\cos\theta e^{\tilde{\tau}_1}+\sin\theta e^{\tilde{\tau}_0}}, \quad \emph{for $e^{i\theta}\in\mathbb{R}P^1$}.
\end{equation}
Note that there might be singularities for some $\tau_{\theta}$ and thus the corresponding Ribaucour sphere congruence may be not regular as the second case in Remark \ref{nonregular}. Generically, the Demoulin family $\{\hat{F}_\theta\}$ can be expressed in terms of the $\mathbb{R}P^1$-family of Ribaucour functions $\{\tau_{\theta}\}$ as
$$\hat{F}_\theta=\mathcal{R}_{\tau_{\theta}}(F).$$
We call such $\{\tau_{\theta}\}$ a \emph{Demoulin family of Ribaucour functions} for $F$.

On the other hand, the condition (\ref{Bianchicond}) also confirms the existence of another two $\nabla$-parallel sections $\hat{\sigma}_i\in\Gamma\mathcal{R}_{s_{\hat{\tau}_i}}(\hat{F}_i)$ $(i=0,1)$, with $\hat{\tau}_i\in\mathcal{R}_{\hat{F}_i}(M)$. Then either of $(\tau_i,\hat{\tau}_i)$, looked as two Ribaucour functions for $\hat{F}_i$, derives the unique dual Demoulin family $\{F_{\theta}=\mathcal{R}_{\hat{\tau}^i_{\theta}}(\hat{F}_i)\}$ through $F$ in the same way as above. Here the uniqueness conclusion is because of the $\nabla$-parallel sections $\hat{\sigma}^i_{\theta}$ consisting of the projective line $[\sigma_i, \hat{\sigma}_i]\subset\mathcal{Q}^2$ all lie on $\Gamma \hat{F}_i$. Now we briefly explain how (\ref{Bianchicond}) gives $\hat{\sigma}_0$ for instance.
Set
$$\hat{\sigma}_0:=v(\hat{\xi}_0-\hat{\tau}_{0}\hat{f}_0-\hat{\tau}_{0}t_0+t_1)=v(\xi-\tau_0f+(\tau_0-\hat{\tau}_0)\hat{f}_0-\hat{\tau}_0t_0+t_1),$$
where the functions $v,\hat{\tau}_{0}$ are to be determined. Then $\hat{\sigma}_0$ is $\nabla$-parallel if and only if
\begin{equation}\label{hatsigma}
(d\hat{\sigma}_0, \hat{f}_1+t_0)=0, \quad (d\hat{\sigma}_0,\sigma_1)=0.
\end{equation}
As in the proof of Theorem $5.4$ of \cite{BH} the orthogonal projection to $V$ of $d\hat{f}_0-(\hat{f}_0+t_0)\hat{\alpha}_{\tau_0}$ takes values in $s_0$ whence we can assume it to be $\gamma(\xi-\tau_0f-\tau_0t_0+t_1)$ for some $1$-form $\gamma$. Then direct computation shows
$$(d\hat{f}_0-(\hat{f}_0+t_0)\hat{\alpha}_{\tau_0}, \hat{f}_1+t_0)=(\tau_1-\tau_0)(a_1-1)\gamma.$$ Using this formula, we get from (\ref{hatsigma}) the following equations for $v,\hat{\tau}_{0}$:
\begin{eqnarray*}
&&-d\ln|v|=d\ln|\tau_0-\hat{\tau}_0|+\hat{\alpha}_{\tau_0},\\
&&d\ln|\tau_0-\hat{\tau}_{0}|=(\alpha_{\tau_1}-\hat{\alpha}_{\tau_0})+d\ln|\tau_1-\tau_0|+(\tau_0-\hat{\tau}_{0})\gamma.
\end{eqnarray*}
Note that $M$ is simply connected and $\alpha_{\tau_1},\hat{\alpha}_{\tau_0}$ are closed by Theorem \ref{Ricond}, to solve them it suffices to prove $$d\Big((\tau_0-\hat{\tau}_{0})\gamma\Big)=0,$$
which can be verified from the formulae above with the help of (\ref{Bianchicond}) and we omit the lengthy computation here.
We conclude these as
\begin{cor}
Given two pointwise distinct Ribaucour functions $\tau_0,\tau_1\in \mathcal{R}_F(M)$ on a simply connected Legedre submanifold satisfying (\ref{Bianchicond}), we have a Demoulin family of Ribaucour functions $\{\tau_{\theta}\}$ for $F$ as in (\ref{Ribaufunctionfamily}) giving the Demoulin family $\{\hat{F}_\theta=\mathcal{R}_{\tau_{\theta}}(F)\}$. The dual Demoulin family $\{F_{\theta}\}$ through $F$ is also determined uniquely.
\end{cor}

At last, we remark that the pair of a constant function $\tau_0$ (as a Ribaucour function corresponding to a parallel transform in the hypersurface case) and any other Ribaucour function $\tau_1\in\mathcal{R}_F(M)$ always satisfies (\ref{Bianchicond}) and thus induces the Demoulin family $\{\hat{F}_{\theta}=\mathcal{R}_{\tau_{\theta}}(F)\}$ of Ribaucour transforms of $F$.
\begin{ack}
 I would like to thank my host Professor Gudlaugur Thorbergsson for his valuable conversations, hospitality and support when I undertake a postdoctoral research fellowship supported by Alexander von Humboldt Foundation in the University of Cologne.
\end{ack}

\end{document}